\documentclass{article}
\usepackage[utf8]{inputenc}
\usepackage{amsthm,amsmath,amssymb,xcolor,placeins}
\usepackage{tikz}
\usepackage{multicol}
\usepackage{enumerate}
\usepackage{xparse}
\usepackage{authblk}
\usepackage{soul}

\usepackage{authblk}
\usepackage{todonotes}

\usepackage{hyperref}
\usetikzlibrary{patterns}

\newtheorem{thm}{Theorem}[section]
\newtheorem{conj}[thm]{Conjecture}
\newtheorem{cor}[thm]{Corollary}

\newtheorem{lemma}[thm]{Lemma}

\theoremstyle{definition}
\newtheorem*{defn}{Definition}

\newtheorem{obs}[thm]{Observation}
\newtheorem{problem}[thm]{Problem}

\NewDocumentEnvironment{manual}{O{obs}m}
 {
  \addtocounter{obs}{-1}
  \renewcommand\theobs{#2}
  \begin{#1}
 }
 {\end{#1}}
 
  \NewDocumentEnvironment{manual2}{O{lemma}m}
 {
  \addtocounter{lemma}{-1}
  \renewcommand\thelemma{#2}
  \begin{#1}
 }
 {\end{#1}}

\renewcommand\det{\operatorname{Det}}
\newcommand\dist{\operatorname{Dist}}
\newcommand\aut{\operatorname{Aut}}
\newcommand\tG{\widetilde G}
\newcommand\tH{\widetilde H}
\newcommand\tT{\widetilde T}
\newcommand\tS{\widetilde S}

\newcommand\tA{\widetilde A}
\newcommand\tB{\widetilde B}

\newcommand\gm{\mu_t(G)}
\newcommand\gh{\mu_t(H)}
\newcommand\mug{\mu(G)}
\newcommand\Tt{T_{t}}

\begin{document}

\title{Symmetry Parameters for Mycielskian Graphs}

\author[1]{Debra Boutin}
\author[2]{Sally Cockburn}
\author[3]{Lauren Keough}
\author[4]{Sarah Loeb}
\author[5]{K.~E. Perry}
\author[6]{Puck Rombach}
\affil[1]{\url{dboutin@hamilton.edu}, Hamilton College, Clinton, NY}
\affil[2]{\url{scockbur@hamilton.edu}, Hamilton College, Clinton, NY}
\affil[3]{\url{keoulaur@gvsu.edu}, Grand Valley State University, Allendale Charter Township, MI}
\affil[4]{\url{sloeb@hsc.edu}, Hampden-Sydney College, Hampden-Sydney, VA}
\affil[5]{\url{kperry@soka.edu}, Soka University of America, Aliso Viejo, CA}
\affil[6]{\url{puck.rombach@uvm.edu}, University of Vermont, Burlington, VT}
\renewcommand\Affilfont{\footnotesize}

\date{\today}

\maketitle

\section{Introduction}\label{sec:intro}

 A coloring of the vertices of a graph $G$ with colors from $\{1,\ldots, d\}$ is called a \emph{$d$-distinguishing coloring} if no nontrivial automorphism of $G$ preserves the color classes. A graph is called \emph{$d$-distinguishable} if it has a $d$-distinguishing coloring. The distinguishing number of $G$, denoted $\dist(G)$, is the smallest number of colors necessary for a distinguishing coloring of $G$. Albertson and Collins introduce graph distinguishing in~\cite{AlCo1996}.  Independently, in~\cite{Ba1977}, Babai defines the same idea, but calls it an \emph{asymmetric coloring}.  Here, we continue the terminology of Albertson and Collins. A substantial amount of work in graph distinguishing in the last few decades proves that, for a large number of graph families, all but a finite number of members are 2-distinguishable. Examples of such families of finite graphs include: hypercubes $Q_n$ with $n\geq 4$~\cite{BoCo2004}, Cartesian powers $G^n$ for a connected graph $G\ne K_2,K_3$ and $n\geq 2$~\cite{Al2005, ImKl2006,KlZh2007}, and Kneser graphs $K_{n:k}$ with $n\geq 6, k\geq 2$~\cite{AlBo2007}. Examples of such families of infinite graphs include: the denumerable random graph~\cite{ImKlTr2007}, the infinite hypercube~\cite{ImKlTr2007}, and denumerable vertex-transitive graphs of connectivity 1~\cite{SmTuWa2012}. 

In 2007, Imrich~\cite{IW} asked if, within 2-distinguishable colorings, we could find the minimum number of times a second color must be used. In response in~\cite{Bo2008}, Boutin defines the \emph{cost of 2-distinguishing} $G$, denoted $\rho(G)$, to be the minimum size of a color class over all 2-distinguishing colorings of $G$. Some of the graph families with known or bounded 2-distinguishing cost are hypercubes with $\lceil \log_2 n \rceil {+} 1 \leq \rho(Q_n) \leq 2\lceil \log_2 n \rceil {-}1$ for $n\geq 5$~\cite{Bo2008}, Kneser graphs with $\rho(K_{2^{m}-1:2^{m{-}1}{-}1}) = m{+}1$~\cite{Bo2013b}, and the Cartesian product of $K_{2^m}$ and an asymmetric graph on $m$ vertices $H$ with $\rho(K_{2^m} \Box H) = m\cdot 2^{m{-}1}$~\cite{BoIm2017}. 
 
A determining set is a useful tool in finding the distinguishing number and, when relevant, the cost of 2-distinguishing. A subset $S$ of $V(G)$ is a \emph{determining set} for $G$ if the only automorphism that fixes the elements of $S$ pointwise is the trivial automorphism. Equivalently, $S$ is a determining set for $G$ if, whenever $\varphi$ and $\psi$ are automorphisms of $G$ with $\varphi(s)=\psi(s)$ for all $s\in S$, then $\varphi=\psi$~\cite{Bo2006}. The \emph{determining number} of a graph $G$, denoted $\det(G)$, is the size of a smallest determining set. Intuitively, if we think of automorphisms of a graph as allowing vertices to move among their relative positions, the determining number is the fewest pins needed to \lq\lq pin down" the graph. 

For some graph families, we only have bounds on the determining number. For instance, for Kneser graphs, $\log_2 (n{+}1)\leq \det(K_{n:k}) \leq n{-}k$, with both bounds sharp~\cite{Bo2006}. However, there are families for which we know the determining numbers of its members exactly. For example, for hypercubes $\det(Q_n){=}\lceil \log_2 n \rceil {+}1$~\cite{Bo2009}, and for generalized Petersen graphs $\det(G(n,k))=2$ if $(n,k)\ne (4,1),(5,2),(10,3)$ and $\det(G(n,k))=3$ otherwise~\cite{Da2020}. 

Though distinguishing numbers and determining numbers were introduced by different people and for different purposes, they have strong connections. Albertson and Boutin show in~\cite{AlBo2007} that if $G$ has a determining set $S$ of size $d$, then giving each vertex in $S$ a distinct color from $\{1,\ldots,d\}$ and coloring the remaining vertices with a $(d+1)^{\text{st}}$ color yields a $(d+1)$-distinguishing coloring of $G$. Thus, $\dist(G) \leq \det(G) {+}1$. Furthermore, for a 2-distinguishing coloring of $G$, Boutin~\cite{Bo2013a} observes that the requirement that only the trivial automorphism preserves the color classes setwise means that only the trivial automorphism preserves them pointwise. Consequently, each of the color classes in a 2-distinguishing coloring is a determining set for the graph, though not necessarily of minimum size. Thus, if $G$ is 2-distinguishable, then $\det(G) \leq \rho(G)$.  In this paper, we refer to the distinguishing number, determining number, and cost of 2-distinguishing, collectively as the \emph{symmetry parameters} of a graph.

In \cite{My1955}, Mycielski introduces a construction that takes a finite simple graph $G$ and produces a larger graph $\mu(G)$ called the (\emph{traditional}) \emph{Mycielskian} of $G$, with the same clique number and a strictly larger chromatic number. In particular, the Myielskian construction can be used to produce families of triangle-free graphs with increasing chromatic number.
A formal definition can be found in Section~\ref{sec:Const}. The \emph{generalized Mycielskian} of graph $G$, denoted $\mu_t(G)$, is defined by Stiebitz in~\cite{St1985} (cited in \cite{Ta2001}) to construct families of graphs with arbitrarily large odd girth and increasing chromatic number.  Independently, Ngoc defines it in~\cite{Va1987} (cited in~\cite{Va1995}). The fundamental difference between the Mycielskian and the generalized Mycielskian of a graph is that the former has 1 level of independent vertices, while the latter has $t$ levels of independent vertices, for some $t\geq 1$.  The definition of $\mu_t(G)$ appears in Section~\ref{sec:Const}.

One strength of the Mycielskian constructions is their ability to build large families of graphs with a given parameter fixed and other parameters strictly growing.  In the last few decades, this has motivated significant work on parameters of $\mu_t(G)$ in terms of the same parameters for $G$. See for example \cite{FiMcBo1998, LWLG2006, ChXi2006, PaZh2010, BaRa2015, AbRa2019, BCKLPR2020a, BCKLPR2020b}.  
In this chapter, for finite simple graphs $G$, we compare symmetry parameters for $G$ and $\mu_t(G)$.

This chapter is organized as follows. Section~\ref{sec:Const} gives the definitions of the Mycielskian constructions. Section~\ref{sec:DistDet} contains a review previous results on symmetry parameters of Mycielskians of graphs, as well as new results on symmetry parameters of the Mycielskian of graphs with isolated vertices. Lastly, Section~\ref{sec:Open} lists some open problems.

\section{Mycielskian Construction}\label{sec:Const}

\begin{figure}[h]
\centering
 \begin{tikzpicture}[scale=.7]

\begin{scope}[shift={(0,0)},scale=.75]
\draw[black!100,line width=1.5pt] (-30:2.5) -- (-150:2.5);
\draw[black!100,line width=1.5pt] (90:2.5) -- (-150:2.5);
\draw[black!100,line width=1.5pt] (-30:2.5) -- (90:2.5);
 \draw[fill=black!100,line width=1] (-30:2.5) circle (.2);
 \draw[fill=black!100,line width=1] (90:2.5) circle (.2);
 \draw[fill=black!100,line width=1] (-150:2.5) circle (.2);
\end{scope}

\begin{scope}[shift={(6,0)},scale=.75]
\draw[black!100,line width=1.5pt] (-30:2.5) -- (-150:2.5);
\draw[black!100,line width=1.5pt] (90:2.5) -- (-150:2.5);
\draw[black!100,line width=1.5pt] (-30:2.5) -- (90:2.5);
\draw[black!100,line width=1.5pt] (-30:2.5) -- (-150:1.25);
\draw[black!100,line width=1.5pt] (90:2.5) -- (-150:1.25);
\draw[black!100,line width=1.5pt] (-30:2.5) -- (90:1.25);
\draw[black!100,line width=1.5pt] (-30:1.25) -- (-150:2.5);
\draw[black!100,line width=1.5pt] (90:1.25) -- (-150:2.5);
\draw[black!100,line width=1.5pt] (-30:1.25) -- (90:2.5);
\draw[black!100,line width=1.5pt] (-30:1.25) -- (0,0);
\draw[black!100,line width=1.5pt] (90:1.25) -- (0,0);
\draw[black!100,line width=1.5pt] (-150:1.25) -- (0,0);
 \draw[fill=black!100,line width=1] (-30:2.5) circle (.2);
 \draw[fill=black!100,line width=1] (90:2.5) circle (.2);
 \draw[fill=black!100,line width=1] (-150:2.5) circle (.2);
 \draw[fill=orange!100,line width=1] (-30:1.25) circle (.2);
 \draw[fill=orange!100,line width=1] (90:1.25) circle (.2);
 \draw[fill=orange!100,line width=1] (-150:1.25) circle (.2);
 \draw[fill=white!100,line width=1] (0,0) circle (.2);
\end{scope}

\begin{scope}[shift={(-1.5,-5)}]
\draw[black!100,line width=1.5pt] (0,0) -- (1,0); 
\draw[fill=black!100,line width=1] (0,0) circle (.15);
\draw[fill=black!100,line width=1] (1,0) circle (.15);
\draw[fill=black!100,line width=1] (2,0) circle (.15);
\draw[fill=black!100,line width=1] (2.5,0) circle (.15);
\draw[fill=black!100,line width=1] (3,0) circle (.15);
\end{scope}
\begin{scope}[shift={(4.5,-5)}]
\draw[black!100,line width=1.5pt] (0,0) -- (1,0); 
\draw[fill=black!100,line width=1] (0,0) circle (.15);
\draw[fill=black!100,line width=1] (1,0) circle (.15);
\draw[fill=black!100,line width=1] (2,0) circle (.15);
\draw[fill=black!100,line width=1] (2.5,0) circle (.15);
\draw[fill=black!100,line width=1] (3,0) circle (.15);
\draw[black!100,line width=1.5pt] (0,0) -- (1,1); 
\draw[black!100,line width=1.5pt] (0,1) -- (1,0); 
\draw[black!100,line width=1.5pt] (1.5,2) -- (0,1); 
\draw[black!100,line width=1.5pt] (1.5,2) -- (1,1); 
\draw[black!100,line width=1.5pt] (1.5,2) -- (2,1); 
\draw[black!100,line width=1.5pt] (1.5,2) -- (2.5,1); 
\draw[black!100,line width=1.5pt] (1.5,2) -- (3,1); 
\draw[fill=orange!100,line width=1] (0,1) circle (.15);
\draw[fill=orange!100,line width=1] (1,1) circle (.15);
\draw[fill=orange!100,line width=1] (2,1) circle (.15);
\draw[fill=orange!100,line width=1] (2.5,1) circle (.15);
\draw[fill=orange!100,line width=1] (3,1) circle (.15);
\draw[fill=white!100,line width=1] (1.5,2) circle (.15);
\end{scope}

\end{tikzpicture}

\caption{Top: $K_3$ and $\mu(K_3)$, drawn with concentric levels with the root in the middle. Bottom: $K_2+3K_1$ and  $\mu(K_2+3K_1)$, drawn with vertical levels with the root at the top.}

\label{fig:GenMyExamplesTrad}

\end{figure}
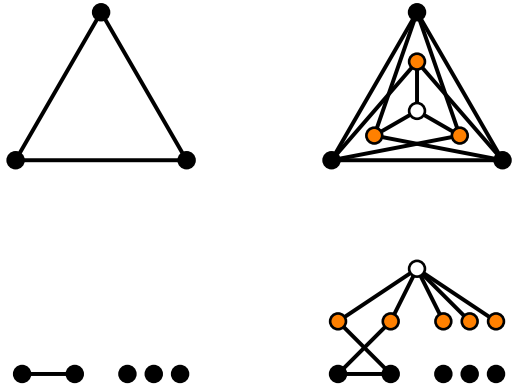

The following is Mycielski's construction for taking a finite simple graph $G$, and producing a larger graph called the \emph{Mycielskian} of $G$. 
\begin{defn}
Let $G$ be a graph with vertex set $V(G) = \{v_1,\dots, v_n\}$. The \emph{Mycielskian of $G$}, denoted $\mu(G)$, has vertex set
\[
V(\mu(G)) = \{ v_1, \dots, v_n, u_1, \dots u_n, w\}.
\]
For each edge $v_iv_j \in E(G)$, $\mu(G)$ has edges $v_iv_j$, $u_iv_j$ and $v_iu_j$; in addition, $u_iw \in E(\mu(G))$ for all $1 \le i \le n$. 
\end{defn}
Note that $\mu(G)$ contains $G$ as an induced subgraph. We refer to the vertices $v_1, \dots , v_n$ in $\mu(G)$ as {\em original vertices} and the vertices $u_1, \dots , u_n$ as {\em shadow vertices}. The vertex $w$ is called the {\em root}.  See Figure~\ref{fig:GenMyExamplesTrad} for illustrations of $\mu(K_3)$ and $\mu(K_2+3K_1)$.

The Mycielskian construction can of course be iterated.  Beginning with $\mu^1(G)=\mu(G)$, define $\mu^k(G)$ inductively as $\mu^k(G)=\mu(\mu^{k-1}(G))$. Mycielski uses this iterated process to define the \emph{classical Mycielskian graphs}, $M_k=\mu^k(K_2)$. Note that indices for this family may differ by publication. This graph family proves the existence of triangle-free graphs with arbitrarily large chromatic number.

The following is the definition of a \emph{generalized Mycielskian} of $G$, also known as a \emph{cone over} $G$.

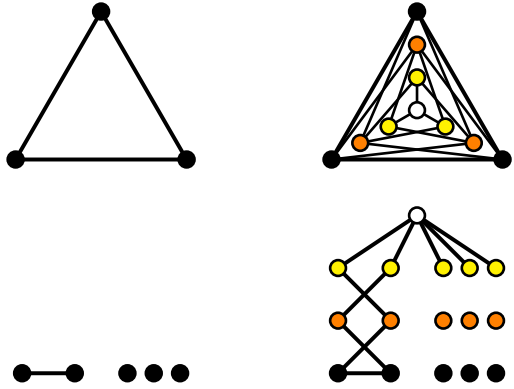
\begin{figure}[h]
\centering
 \begin{tikzpicture}[scale=.7]

\begin{scope}[shift={(0,0)},scale=.75]
\draw[black!100,line width=1.5pt] (-30:2.5) -- (-150:2.5);
\draw[black!100,line width=1.5pt] (90:2.5) -- (-150:2.5);
\draw[black!100,line width=1.5pt] (-30:2.5) -- (90:2.5);
 \draw[fill=black!100,line width=1] (-30:2.5) circle (.2);
 \draw[fill=black!100,line width=1] (90:2.5) circle (.2);
 \draw[fill=black!100,line width=1] (-150:2.5) circle (.2);
\end{scope}

\begin{scope}[shift={(6,0)},scale=.75]
\draw[black!100,line width=1.5pt] (-30:2.5) -- (-150:2.5);
\draw[black!100,line width=1.5pt] (90:2.5) -- (-150:2.5);
\draw[black!100,line width=1.5pt] (-30:2.5) -- (90:2.5);
\draw[black!100,line width=1pt] (-30:2.5) -- (-150:1.66);
\draw[black!100,line width=1pt] (90:2.5) -- (-150:1.66);
\draw[black!100,line width=1pt] (-30:2.5) -- (90:1.66);
\draw[black!100,line width=1pt] (-30:1.66) -- (-150:2.5);
\draw[black!100,line width=1pt] (90:1.66) -- (-150:2.5);
\draw[black!100,line width=1pt] (-30:1.66) -- (90:2.5);
\draw[black!100,line width=1pt] (-30:.83) -- (-150:1.66);
\draw[black!100,line width=1pt] (90:.83) -- (-150:1.66);
\draw[black!100,line width=1pt] (-30:.83) -- (90:1.66);
\draw[black!100,line width=1pt] (-30:1.66) -- (-150:.83);
\draw[black!100,line width=1pt] (90:1.66) -- (-150:.83);
\draw[black!100,line width=1pt] (-30:1.66) -- (90:.83);
\draw[black!100,line width=1pt] (-30:.83) -- (0,0);
\draw[black!100,line width=1pt] (90:.83) -- (0,0);
\draw[black!100,line width=1pt] (-150:.83) -- (0,0);
 \draw[fill=black!100,line width=1] (-30:2.5) circle (.2);
 \draw[fill=black!100,line width=1] (90:2.5) circle (.2);
 \draw[fill=black!100,line width=1] (-150:2.5) circle (.2);
 \draw[fill=yellow!100,line width=1] (-30:.83) circle (.2);
 \draw[fill=yellow!100,line width=1] (90:.83) circle (.2);
 \draw[fill=yellow!100,line width=1] (-150:.83) circle (.2);
  \draw[fill=orange!100,line width=1] (-30:1.66) circle (.2);
 \draw[fill=orange!100,line width=1] (90:1.66) circle (.2);
 \draw[fill=orange!100,line width=1] (-150:1.66) circle (.2);
 \draw[fill=white!100,line width=1] (0,0) circle (.2);
\end{scope}

\begin{scope}[shift={(-1.5,-5)}]
\draw[black!100,line width=1.5pt] (0,0) -- (1,0); 
\draw[fill=black!100,line width=1] (0,0) circle (.15);
\draw[fill=black!100,line width=1] (1,0) circle (.15);
\draw[fill=black!100,line width=1] (2,0) circle (.15);
\draw[fill=black!100,line width=1] (2.5,0) circle (.15);
\draw[fill=black!100,line width=1] (3,0) circle (.15);
\end{scope}
\begin{scope}[shift={(4.5,-5)}]
\draw[black!100,line width=1.5pt] (0,0) -- (1,0); 
\draw[fill=black!100,line width=1] (0,0) circle (.15);
\draw[fill=black!100,line width=1] (1,0) circle (.15);
\draw[fill=black!100,line width=1] (2,0) circle (.15);
\draw[fill=black!100,line width=1] (2.5,0) circle (.15);
\draw[fill=black!100,line width=1] (3,0) circle (.15);
\draw[black!100,line width=1.5pt] (0,0) -- (1,1); 
\draw[black!100,line width=1.5pt] (0,1) -- (1,0);
\draw[black!100,line width=1.5pt] (0,1) -- (1,2); 
\draw[black!100,line width=1.5pt] (0,2) -- (1,1);
\draw[black!100,line width=1.5pt] (1.5,3) -- (0,2); 
\draw[black!100,line width=1.5pt] (1.5,3) -- (1,2); 
\draw[black!100,line width=1.5pt] (1.5,3) -- (2,2); 
\draw[black!100,line width=1.5pt] (1.5,3) -- (2.5,2); 
\draw[black!100,line width=1.5pt] (1.5,3) -- (3,2);
\draw[fill=orange!100,line width=1] (0,1) circle (.15);
\draw[fill=orange!100,line width=1] (1,1) circle (.15);
\draw[fill=orange!100,line width=1] (2,1) circle (.15);
\draw[fill=orange!100,line width=1] (2.5,1) circle (.15);
\draw[fill=orange!100,line width=1] (3,1) circle (.15);
\draw[fill=yellow!100,line width=1] (0,2) circle (.15);
\draw[fill=yellow!100,line width=1] (1,2) circle (.15);
\draw[fill=yellow!100,line width=1] (2,2) circle (.15);
\draw[fill=yellow!100,line width=1] (2.5,2) circle (.15);
\draw[fill=yellow!100,line width=1] (3,2) circle (.15);
\draw[fill=white!100,line width=1] (1.5,3) circle (.15);
\end{scope}

\end{tikzpicture}

\caption{Top: $K_3$, $\mu_2(K_3)$, drawn with concentric levels with the root in the middle. Bottom: $K_2+3K_1$ and $\mu_2(K_2+3K_1)$, drawn with vertical levels with the root at the top. }
\label{fig:GenMyExamplesGen}

\end{figure}

\begin{defn}
 Let $G$ be a graph with vertex set $V(G)=\{v_1,\dots,v_n\}$ and let $t \geq 1$. The \emph{generalized Mycielskian of the graph}, denoted $\gm$, has vertex set
\[V(\gm)=\{u^0_1,\ldots,u^0_n,u^1_1,\ldots,u^1_n,\dots,u^t_1,\ldots,u^t_n,w\}.\]
 For each edge $v_i v_j$ in $G$, the graph $\mu_t(G)$ has edge $v_i v_j  = u_i^0u_j^0$, as well as edges $u^s_iu^{s{+}1}_j$ and $u^s_j u^{s+1}_i$, for $0\leq s <t$. Finally, $\mu_t(G)$ has edges $u^t_i w$ for $1 \leq i \leq n$. 
\end{defn}
We say that vertex $u_i^s$ is \emph{at level $s$}. In addition, we make the identification $u_i^0 = v_i$ and refer to the vertices at level zero as \emph{original vertices} and the vertices at level $s \ge 1$ as \emph{shadow vertices} and to $w$ as the \emph{root}. Notice that the root is only adjacent to the shadow vertices at level $t$. We call level $t$ the \emph{top level}. Observe that $\mu_1(G) = \mug$. Thus, for ease of notation, we omit the subscript when $t=1$.  See Figure~\ref{fig:GenMyExamplesGen} for illustrations of $\mu_2(K_3)$ and $\mu_2(K_2+3K_1)$.
 
 To help understand the structure of the Mycielskian and generalized Mycielskian, it is useful to understand the degrees of the vertices. Denote the degree of vertex $v\in V(G)$ by $\deg_G(v)$ and the degree of $x\in V(\mu(G))$ by $\deg_{\mu(G)}(x)$. By the generalized Mycielskian construction, for each original vertex $u\in V(G)$, $\deg_{\mu_t(G)}(u) = 2\deg_G(u)$.  Further, for $1\leq i \leq t-1$, for the shadow of $u$ at level $i$, $\deg_{\mu_t(G)}(u^i) = 2\deg_G(u)$. However, for the shadow of $u$ at the top level $t$, $\deg_{\mu_t(G)}(u^t) = \deg_G(u) +1$.  Finally, $\deg_{\mu_t(G)}(w)=|V(G)|$.

To address symmetry parameters of $\mu_t(G)$, it is useful to be able to discuss the automorphisms of $\mu_t(G)$ in terms of the automorphisms of $G$.  Boutin, Cockburn, Keough, Loeb, Perry and Rombach prove in \cite{BCKLPR2020a, BCKLPR2020b} that if $G$ is not a star graph, that is, not $K_{1,m}$ for any $m\geq 0$, then each automorphism of $\mu_t(G)$ fixes the root, preserves the levels of the vertices, and restricts to an automorphism of $G$.  This is stated formally in the following lemma.

\begin{lemma} \label{lem:phiwisv} {\rm \cite{BCKLPR2020a, BCKLPR2020b}}  Let 
$G \neq K_{1,m}$ for any $m \geq 0$, and suppose that $G$ has $n$ vertices and no isolated vertices.  Let  $t \geq 1$, and  let $\widehat \alpha$ be an automorphism of $\mu_t(G)$. Then, 

\begin{enumerate}[(i)]

\item $\widehat\alpha(w) = w$,

\item $\widehat \alpha$ preserves the level of vertices, that is,    $\widehat\alpha(\{u_1^s \dots , u_n^s\}) \subseteq \{u_1^s \dots , u_n^s\}$ for all $0 \le s \le t$, and
 
\item $\widehat\alpha$ restricted to $\{u_1^0, \dots , u_n^0\} = \{v_1, \dots, v_n\}$ is an automorphism of $G$.
\end{enumerate}
\end{lemma}

\section{Distinguishing and Determining Mycielskians} \label{sec:DistDet}

In ~\cite{AlSo2019b}, Alikhani and Soltani show that the classical Mycielskian graphs defined in Section~\ref{sec:Const} satisfy $\dist(M_k) = 2$ for all $k \ge 2$. To find symmetry parameters of Mycielskians of arbitrary graphs, they consider the role of twin vertices. Two vertices in a graph are said to be \emph{twins} if they have the same open neighborhood, and a graph is said to be \emph{twin-free} if it does not contain any twins. In particular, Alikhani and Soltani prove  that if $G$ is twin-free with at least two vertices, then $\dist(\mu(G)) \le \dist(G){+}1$.  Further, in~\cite{AlSo2019b} the authors conjecture that for all but a finite number of connected graphs $G$ with at least 3 vertices, $\dist(\mu(G)) \le \dist(G)$. In \cite{BCKLPR2020a}, Boutin, Cockburn, Keough, Loeb, Perry, and Rombach prove the Alikhani and Soltani conjecture, and follow up in \cite{BCKLPR2020b} by studying the determining number and cost of $2$-distinguishing for $\mu(G)$ and $\mu^t(G)$.  The main results are summarized in Theorems~\ref{thm:distmut}-\ref{thm:DetMain}.

\begin{thm}\label{thm:distmut} {\rm \cite{BCKLPR2020a}} Let $G\neq K_1, K_2$  be a graph with $\ell \geq 0$  isolated vertices. If $t \ell > \dist(G)$, then $\dist(\mu_t(G)) = t \ell$; otherwise, $\dist(\mu_t(G))\leq \dist(G)$. \end{thm}
  
\begin{thm}\label{thm:twinfreeDet1}{\rm \cite{BCKLPR2020b}}
 Let $G$ be a twin-free graph with no isolated vertices such that $\det(G) \geq 2$. Then for $t \geq \det(G)-1 $,
 
 \begin{itemize}  
\item[(i)] $\dist(\gm) = 2$,   \item[(ii)] $\det(\gm) = \det(G)$, and \item[(iii)] $\rho(\gm) = \det(G)$.
\end{itemize}
\end{thm}

When $G$ has twins, but no isolated vertices, the same authors characterize the determining number of $\gm$. To understand this result, we define a \emph{minimum twin cover} $T$ of $G$ as a set consisting of all but one vertex from each set of mutually twin vertices. See Section~\ref{sec:isowtwins} for more details.

\begin{thm}\label{thm:DetMain}{\rm \cite{BCKLPR2020b}}
 Let $G$ be a graph with no isolated vertices. Let $T$ be a (possibly empty) minimum twin cover of $G$. Let $t \geq 1$.
  
 \begin{enumerate}[(i)]
  
   \item If $G = K_2$, then $\det(G) = 1$ and $\det(\gm) = 2.$ 

    \item If $G \neq K_2$, then $\det(\gm) = t|T| + \det(G).$
    \end{enumerate}

\end{thm}

\begin{figure}[h]
    \centering
\begin{tikzpicture}[scale=.7]

\draw[black!100,line width=1.5pt] (-1,0) -- (1,0); 
\draw[black!100,line width=1.5pt] (-1,0) to[out=-40,in=-140]  (1,0); 
\draw[fill=black!100,line width=1] (-1,0) circle (.15);
\draw[fill=black!100,line width=1] (0,0) circle (.15);
\draw[fill=black!100,line width=1] (1,0) circle (.15);
\draw[fill=black!100,line width=1] (2,0) circle (.15);
\draw[fill=black!100,line width=1] (2.5,0) circle (.15);
\draw[fill=black!100,line width=1] (3,0) circle (.15);

\begin{scope}[shift={(6,0)}]
\draw[black!100,line width=1.5pt] (-1,0) -- (1,0); 
\draw[black!100,line width=1.5pt] (-1,0) to[out=-40,in=-140]  (1,0); 
\draw[fill=black!100,line width=1] (-1,0) circle (.15);
\draw[fill=black!100,line width=1] (0,0) circle (.15);
\draw[fill=black!100,line width=1] (1,0) circle (.15);
\draw[fill=black!100,line width=1] (2,0) circle (.15);
\draw[fill=black!100,line width=1] (2.5,0) circle (.15);
\draw[fill=black!100,line width=1] (3,0) circle (.15);
\draw[black!100,line width=1.5pt] (0,0) -- (1,1); 
\draw[black!100,line width=1.5pt] (0,1) -- (1,0); 
\draw[black!100,line width=1.5pt] (0,0) -- (-1,1); 
\draw[black!100,line width=1.5pt] (1,1) -- (-1,0); 
\draw[black!100,line width=1.5pt] (0,1) -- (-1,0);
\draw[black!100,line width=1.5pt] (-1,1) -- (1,0);
\draw[black!100,line width=1.5pt] (1.5,2) -- (-1,1);
\draw[black!100,line width=1.5pt] (1.5,2) -- (0,1); 
\draw[black!100,line width=1.5pt] (1.5,2) -- (1,1); 
\draw[black!100,line width=1.5pt] (1.5,2) -- (2,1); 
\draw[black!100,line width=1.5pt] (1.5,2) -- (2.5,1); 
\draw[black!100,line width=1.5pt] (1.5,2) -- (3,1); 
\draw[fill=orange!100,line width=1] (-1,1) circle (.15);
\draw[fill=orange!100,line width=1] (0,1) circle (.15);
\draw[fill=orange!100,line width=1] (1,1) circle (.15);
\draw[fill=orange!100,line width=1] (2,1) circle (.15);
\draw[fill=orange!100,line width=1] (2.5,1) circle (.15);
\draw[fill=orange!100,line width=1] (3,1) circle (.15);
\draw[fill=white!100,line width=1] (1.5,2) circle (.15);
\end{scope}
\begin{scope}[shift={(12,0)}]
\draw[black!100,line width=1.5pt] (-1,0) -- (1,0); 
\draw[black!100,line width=1.5pt] (-1,0) to[out=-40,in=-140]  (1,0); 
\draw[fill=black!100,line width=1] (-1,0) circle (.15);
\draw[fill=black!100,line width=1] (0,0) circle (.15);
\draw[fill=black!100,line width=1] (1,0) circle (.15);
\draw[fill=black!100,line width=1] (2,0) circle (.15);
\draw[fill=black!100,line width=1] (2.5,0) circle (.15);
\draw[fill=black!100,line width=1] (3,0) circle (.15);
\draw[black!100,line width=1.5pt] (0,0) -- (1,1); 
\draw[black!100,line width=1.5pt] (0,1) -- (1,0); 
\draw[black!100,line width=1.5pt] (0,0) -- (-1,1); 
\draw[black!100,line width=1.5pt] (1,1) -- (-1,0); 
\draw[black!100,line width=1.5pt] (0,1) -- (-1,0);
\draw[black!100,line width=1.5pt] (-1,1) -- (1,0);

\draw[black!100,line width=1.5pt] (0,1) -- (1,2); 
\draw[black!100,line width=1.5pt] (0,2) -- (1,1); 
\draw[black!100,line width=1.5pt] (0,1) -- (-1,2); 
\draw[black!100,line width=1.5pt] (1,2) -- (-1,1); 
\draw[black!100,line width=1.5pt] (0,2) -- (-1,1);
\draw[black!100,line width=1.5pt] (-1,2) -- (1,1);
\draw[black!100,line width=1.5pt] (1.5,3) -- (-1,2);
\draw[black!100,line width=1.5pt] (1.5,3) -- (0,2); 
\draw[black!100,line width=1.5pt] (1.5,3) -- (1,2); 
\draw[black!100,line width=1.5pt] (1.5,3) -- (2,2); 
\draw[black!100,line width=1.5pt] (1.5,3) -- (2.5,2); 
\draw[black!100,line width=1.5pt] (1.5,3) -- (3,2); 
\draw[fill=orange!100,line width=1] (-1,1) circle (.15);
\draw[fill=orange!100,line width=1] (0,1) circle (.15);
\draw[fill=orange!100,line width=1] (1,1) circle (.15);
\draw[fill=orange!100,line width=1] (2,1) circle (.15);
\draw[fill=orange!100,line width=1] (2.5,1) circle (.15);
\draw[fill=orange!100,line width=1] (3,1) circle (.15);
\draw[fill=yellow!100,line width=1] (-1,2) circle (.15);
\draw[fill=yellow!100,line width=1] (0,2) circle (.15);
\draw[fill=yellow!100,line width=1] (1,2) circle (.15);
\draw[fill=yellow!100,line width=1] (2,2) circle (.15);
\draw[fill=yellow!100,line width=1] (2.5,2) circle (.15);
\draw[fill=yellow!100,line width=1] (3,2) circle (.15);
\draw[fill=white!100,line width=1] (1.5,3) circle (.15);
\end{scope}
\end{tikzpicture}
    \caption{The graphs $G=K_3+3K_1$, $\mu (G)$ and $\mu_2 (G)$. }
    \label{fig:K2tK1ex1}
\end{figure}
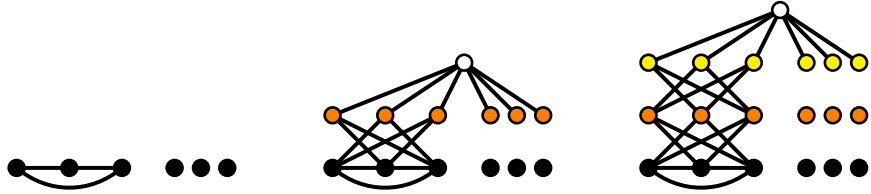

As seen in the results of \cite{BCKLPR2020a, BCKLPR2020b, FiMcBo1998}, isolated vertices can play an important role when investigating parameters of Mycielskians and generalized Mycielskians. If a graph $G$ has isolated vertices, then $\mu_t(G)$ has one component that consists of the generalized Mycielskian of the graph induced by the non-isolated vertices of $G$, and the top-level shadows of isolated vertices as pendant vertices adjacent to the root. However, $\mu_t(G)$ also has $t$ isolated vertices for every isolated vertex in $G$. More formally, if $G = H + \ell K_1$, where $H$ is a graph with at least one edge and no isolated vertices, then $\gm = C + t\ell K_1$, where $C$ is a connected component consisting of $\gh$ with $\ell$ additional pendant vertices adjacent to $w$. See Figure~\ref{fig:K2tK1ex1} which shows the generalized Mycielskians of $G = K_3 + 3K_1$ when $t = 1$ and $2$.

We next provide a lemma that helps us understand the automorphisms and determining numbers of graphs with isolated vertices. We then investigate the role that isolated vertices play when considering the distinguishing and determining number of graphs with and without twin vertices. 
 
\begin{lemma}~\label{lem:detghvsC}
Let $H$ be a graph with at least one edge and no isolated vertices. For $t \geq 1$, let $C$ be a copy of $\gh$ with an additional $\ell \ge 1$ pendant vertices adjacent to $w$: $x_1, x_2, \dots, x_{\ell}$. 
 If $S$ is a minimum size subset of $V(\gh) \setminus \{w\}$ such that $S
\cup \{w\}$ is a determining set for $\gh$, then $S \cup \{x_2, \dots, x_{\ell} \}$ is a minimum size determining set for $C$. 

\end{lemma}

\begin{proof}
Let $S$ be a subset of  $V(\gh) \setminus\{w\}$  of minimum size such that  $S \cup \{ w \}$ is a determining set for $\gh$. 
Observe that $C$ has exactly $\ell$ vertices of degree $1$ and $w$ is the only vertex adjacent to all of these vertices. Thus, $w$ is fixed by every automorphism of $C.$
 If $\ell = 1$, then since $C$ has only one pendant vertex, it must be fixed under every automorphism.
If $\ell \ge 2$, then any automorphism of $C$ that fixes $\{x_2,\dots, x_\ell\}$ must also fix $x_1$. In both cases, any $\alpha \in \aut(C)$ that fixes $\{x_2,\dots, x_\ell\}$ can act nontrivially on only vertices in $V(\gh) \setminus\{w\}$. If $\alpha$ also fixes the vertices in  $S$, then by the assumption that $S \cup \{w\}$ is a determining set for $\gh$, $\alpha$ must be the identity.
Hence, $S\cup\{x_2,\dots, x_\ell\}$ is a determining set for $C$.

To show minimality, suppose that $B$ is a minimum size determining set for $C$ such that $|B| < |S\cup\{x_2,\dots, x_\ell\}|$. Since $x_1,\dots, x_\ell$ are mutually twin, we can assume without loss of generality that $\{x_2,\dots, x_\ell\}\subseteq B$. As already noted, any automorphism of $C$ fixing these vertices also fixes $w$ and $x_1$, and so by the assumption of minimality, $w, x_1 \notin B$. Then $B' = B \setminus \{x_2,\dots, x_\ell\}$ is a subset of $V(\gh) \setminus\{w\}$ and $|B'| < |S|$. 
Now, let $\alpha$ be any automorphism of $\gh$ fixing $B'\cup \{w\}$. Since $\alpha$ can be extended to an automorphism of $C$ fixing $B$, a determining set for $C$, it follows that $\alpha$ restricted to $\gh$ must be the trivial automorphism and so $B'\cup \{w\}$ is a determining set for $\gh$, a contradiction.
\end{proof}

 Notice that Lemma~\ref{lem:detghvsC} shows that if $G = H+ K_1$ and $w$ is not in any minimum determining set for $\gh$, then $\det(\gh) = \det(C)$. By Lemma~\ref{lem:phiwisv}, if $H \neq K_{1,m}$ for any $m\geq 0$, then all automorphisms of $\gh$ fix $w$. Thus, if $H$ is not a star graph, then $\det(\gh) = \det(C)$.  

\subsection{\boldmath An Isolated Vertex in Twin-Free $G$} \label{sec:isonotwins}

Recall that any collection of vertices are said to be twins if they all share the same open neighborhood. Since there are automorphisms that permute twin vertices while leaving other vertices fixed, the presence of twin vertices affects the determining and distinguishing numbers. Thus, we address graphs with twin vertices separately in Section~\ref{sec:isowtwins}.  In this section, we extend results from both~\cite{BCKLPR2020a} and~\cite{BCKLPR2020b} on the determining and  distinguishing numbers of twin-free graphs to twin-free graphs with an isolated vertex. Note that a twin-free graph can have at most one isolated vertex. 

In~\cite{BCKLPR2020a}, Boutin, Cockburn, Keough, Loeb, Perry and Rombach show that $\dist(\gm) \leq \dist(G)$ for most graphs $G$. However, in~\cite{BCKLPR2020b} the same authors show that for most graphs equality holds for the determining number of twin-free graphs with no isolated vertices. In particular, they show the following.

\begin{thm} \label{twinfreeDet}~{\rm \cite{BCKLPR2020b}}
 Let $G$ be a twin-free graph with no isolated vertices and let $t \geq 1$.
\begin{enumerate}[(i)]
    
     \item If $G = K_2$ then, $\det(G) = 1$ and $\det(\gm) = 2.$ 
     
    \item If $G \neq K_2$, then any minimum size determining set for $G$ is a minimum size determining set for $\gm$ and \[\det(\gm) = \det(G).\]
\end{enumerate}
\end{thm}

The theorem below gives analogous results for  twin-free graphs having exactly one isolated vertex. 

\begin{thm}\label{twinfreeDetIso}
 Let $G$ be a twin-free graph and let $t \geq 1$.
\begin{enumerate}[(i)]
   \item If $G = K_1$, then $\det(G) = 0$ and $\det(\gm) = t.$
    \item If $G$ is of the form $H + K_1$ for some graph $H$ with at least one edge, then \[\det(\gm) = \det(G) + t - 1.\]
\end{enumerate}
\end{thm}

\begin{proof}
If $G = K_1$, then $\det(G) = 0$ since the only automorphism of $G$ is the trivial one.  Note that $\mu_t(K_1) = K_2 + tK_1$ 
and a minimum size determining set consists of one vertex in  $K_2$ and $t-1$ of the isolated vertices. 

Next, suppose $G \neq K_1$ is of the form $H + K_1$, where $H$ does not have an isolated vertex because $G$ is twin-free. Then $\det(G) = \det(H)$.  Moreover, $\gm = C+ tK_1$, where $C$ is a connected component consisting of $\mu_t(H)$ with an additional pendant vertex adjacent to $w$. Thus,
\[
\det(\gm) = \det(C) + t-1.
\]
It therefore suffices to show that $\det(C) = \det(H)$. We divide into two cases.

 If $H = K_2$, then $\mu_t(K_2) = C_{2t+3}$. As noted earlier, any two vertices on an odd cycle constitute a minimum size determining set. Hence if $y\ne w$ is any vertex on the cycle, we can apply Lemma~\ref{lem:detghvsC} to conclude that $S = \{y\}$ is a minimize size determining set for $C$, which implies $\det(C) = 1 = \det(K_2)$.

Now, suppose $H \neq K_2$. Since in addition $G$ is  twin-free,  $H \neq K_{1,m}$ for any $m \ge 0$. Since $H$ is not a star graph, by Lemma~\ref{lem:phiwisv}, every automorphism of $\gh$ fixes $w$, which means that for every $S \subseteq V(\gh) \setminus \{w\}$, $S \cup \{w\}$ is a determining set for $\gh$ only if $S$ is. By Lemma~\ref{lem:detghvsC}, then, if $S$ is a minimum size determining set for $\gh$, it is also one for $C$ and so $\det(C) = \det(\mu_t(H))$.
Applying Theorem~\ref{twinfreeDet}(ii) to $H$ gives $\det(\mu_t(H)) = \det(H)$ and we are done.\end{proof}

Theorem~\ref{twinfreeDet} and Theorem~\ref{twinfreeDetIso} effectively tell us that if $G \neq K_1$ or $K_2$ and $G$ is twin-free, then only two behaviors are possible.

\begin{cor}\label{cor:DetGMtwinfree}
If $G\neq K_1, K_2$ is a twin-free graph, then for all $t \geq 1$, 
\[
\det(\gm) = \begin{cases}
\det(G) & \text{if $G$ does not have an isolated vertex;}\\
 \det(G) + t -1 &\text{otherwise}
 \end{cases}
.\]
\end{cor}

The theorem below, from~\cite{BCKLPR2020b}, addresses the distinguishing number of Mycielskians of twin-free graphs with no isolated vertices. It asserts that for $t$ sufficiently large,  these graphs are $2$-distinguishable and provides information on the cost of $2$-distinguishing. Note that, if $\det(G)=1$, then $\dist(\gm) = 2$ and $\rho(\gm) = 1$.

\begin{thm}~{\rm \cite{BCKLPR2020b}}
Let $G$ be a twin-free graph with no isolated vertices. If  $\det(G) = k \ge 2$, then the following hold.
\begin{enumerate}[(i)]
    \item If $ \lceil \log_2(k+1)\rceil -1 \le t $, then  $\dist(\gm) = 2$ and
\[\rho(\gm) \leq \frac{(k+1) \lceil \log_2(k+1) \rceil}{2}.\]
\item If  $t\ge k-1$, then  $\dist(\gm) = 2$ and $\rho(\gm) = k$.
\end{enumerate}
In (i), the bound on $t$ and on $\rho(\gm)$ are both sharp.
\end{thm}

If a twin-free graph $G$ has one isolated vertex, then $\gm$ has $t$ isolated vertices, which means that for $t \ge 2$, $\gm$ is no longer twin-free. As shown in Theorem~\ref{twinfreeDetIso}, this causes the determining number to grow linearly with $t$.  In any distinguishing coloring, mutually twin vertices must receive different colors and so it is not surprising that the distinguishing number also grows linearly with $t$.

\begin{thm}
Let $G$ be a twin-free graph of the form $H + K_1$ with $\det(G) = k \geq 1$. If $ t \geq \lceil \log_2(k+1)\rceil -1$, then $\dist(\gm) = \max(2,t)$. If in addition $t=1$ or $2$, then $\rho(\gm) =  k + t-1$. 
\end{thm}

\begin{proof}
As noted in the proof of Theorem~\ref{twinfreeDetIso}, $\gm = C + tK_1$ where $C$ consists of $\mu_t(H)$ with one pendant vertex adjacent to $w$. 
Since $H$ is twin-free with no isolated vertices, we can apply  Theorem~\ref{thm:twinfreeDet1} to conclude that $\dist (\mu_t(H)) = 2$. 
We can extend any 2-distinguishing coloring of $\mu_t(H)$ to a 2-coloring of $C$ by coloring the pendant vertex either of the two colors. Every automorphism of $C$ must fix the pendant vertex; every automorphism $\alpha$ of $C$ that fixes these 2 colors classes of $C$ therefore restricts to an automorphism of $\mu_t(H)$ that fixes the 2 colors classes in the 2-distinguishing coloring. By definition, the restriction of $\alpha$ acts as the identity on $\mu_t(H)$ and therefore $\alpha$ is the identity on $C$. Hence this 2-coloring is distinguishing on $C$ and so $\dist(C) \le 2$. From the proof of Theorem~\ref{twinfreeDetIso}, $\det(C) = \det(H) = \det(G) \ge 1$. This implies that $C$ has nontrivial automorphisms and so $\dist (C) \ge 2$. Thus, $\dist(C) = 2$.

A distinguishing coloring for $\gm = C +tK_1$ must assign a different color to each isolated vertex. By reusing two colors used in a  2-distinguishing coloring of $C$, we need at most $\max(0,t-2)$ additional colors.

If $t=1$ or $2$, then $\gm$ is 2-distinguishable and we can compute the cost. By Theorem~\ref{thm:twinfreeDet1} and the fact that $\det(H) = \det(G) = k$, we know $\rho(\mu_t(H)) = k$. Let $S \subset V(\mu_t(H))$ be a set of $k$ vertices in a smallest color class of a 2-distinguishing coloring, and let $c_1$ be the color used for $S$, with $c_2$  being the other color. Since the pendant vertex in $C$ is fixed by any automorphism of $C$, we can color it with  $c_2$.
Thus, $\rho(C) = \rho(\mu_t(H)) = k$. If $t=1$, then we can also color the one isolated vertex in $\mu_1(G) = \mu(G) = C + K_1$ with $c_2$ and so $\rho(\gm) =  k.$ If $t=2$, then we must use both colors on the isolated vertices and so $\rho(\gm) = k +1$.\end{proof}

\subsection{\boldmath Isolated Vertices in $G$ with Twins} \label{sec:isowtwins}

 For an analysis of determining numbers of Mycielskians of graphs with twin vertices and without isolated vertices, see~\cite{BCKLPR2020b}. Here we consider graphs with twin vertices that may have isolated vertices. The proofs in this section are similar to~\cite{BCKLPR2020b}, and we begin with a more condensed version of the background here. 
 
  The following useful observation about the relationship between twin vertices in $G$ and those is  $\gm$ appears in~\cite{BCKLPR2020b}.
 \begin{obs}\label{obs:twins} {\rm \cite{BCKLPR2020b}} \it If original vertices $v_i$ and $v_j$ are twins in $G$, then $u_i^s$ and $u_j^s$ are twins in  $\mu_t(G)$ for all $0 \le s \le t$. Similarly, if $\{u_i^s,u_j^s\}$ are twins in $\mu_t(G)$ for some $0 \le s \le t$, then they are twins for all $0\le s \le t$; in particular, $\{v_i, v_j\}$ are twins in $G$.
 \end{obs}
  
 Let $x \sim y$ be the equivalence relation on $V(G)$ that indicates that $x$ and $y$ are twin vertices. The quotient graph for relation $\sim$, denoted $\tG$, has as its vertices the set of equivalence classes of the form $[x] = \{y \in V(G) \mid x \sim y\}$ with $[x]$ adjacent to $[z]$ in $\tG$ if and only if there exist $p \in [x]$ and $q \in [z]$ such that $pq \in E(G)$. We see that
  \[
  N_{\tG} ([x]) = \{ [z] \mid z \in N_G(x)\}.
  \]
Furthermore, note that $\tG$ is twin-free, and, if $G$ is twin-free, then $\tG=G$. 
  
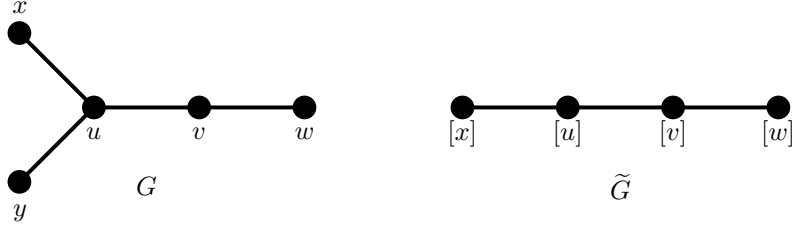
\begin{figure}
    \centering
         \begin{tikzpicture}[scale=.7]
\draw[fill=black!100,line width=1] (0,0) circle (.2);
\draw[fill=black!100,line width=1] (135:2) circle (.2);
\draw[fill=black!100,line width=1] (-135:2) circle (.2);
\draw[fill=black!100,line width=1] (2,0) circle (.2);
\draw[fill=black!100,line width=1] (4,0) circle (.2);
\draw[black!100,line width=1.5pt] (0,0) -- (4,0);
\draw[black!100,line width=1.5pt] (0,0) -- (135:2);
\draw[black!100,line width=1.5pt] (0,0) -- (-135:2);
\draw (0,-.5) node{$u$};
\draw (2,-.5) node{$v$};
\draw (4,-.5) node{$w$};
\draw (-1.4,-2) node{$y$};
\draw (-1.4,1.9) node{$x$};
\draw (1,-1.5) node{$G$};
\begin{scope}[shift={(9,0)}]
\draw[fill=black!100,line width=1] (0,0) circle (.2);
\draw[fill=black!100,line width=1] (-2,0) circle (.2);
\draw[fill=black!100,line width=1] (2,0) circle (.2);
\draw[fill=black!100,line width=1] (4,0) circle (.2);
\draw[black!100,line width=1.5pt] (-2,0) -- (4,0);
\draw (0,-.5) node{$[u]$};
\draw (2,-.5) node{$[v]$};
\draw (4,-.5) node{$[w]$};
\draw (-2,-.5) node{$[x]$};
\draw (1,-1.5) node{$\tG$};
\end{scope}
\end{tikzpicture}
    \caption{A graph $G$ and its quotient graph $\tG$~\cite{BCKLPR2020b}.} 
    \label{fig:tGisP4}
\end{figure}

Since automorphisms preserve neighborhoods, every automorphism $\alpha$ of $G$ induces an automorphism $\widetilde \alpha$ of $\tG$ given by $\widetilde \alpha([x]) = [\alpha(x)]$. However, there may be automorphisms of $\tG$ that are not of this form, as is shown in Figure~\ref{fig:tGisP4}~\cite{BCKLPR2020b} where the only nontrivial automorphism of $G$ is the one interchanging the twin vertices $x$ and $y$, but $\tG = P_4$ has a reflectional nontrivial automorphism.

 Throughout the rest of this section, we use $\tG$ and $\tH$ to denote the quotient graphs of graphs $G$ and $H$, respectively, and $\widetilde \alpha$ to denote an automorphism of a quotient graph. Vertex sets with tilde notation represent subsets of vertices in a quotient graph. In particular, if $B \subseteq V(G)$, then $\tB = \{ [x] \in V(\tG) \mid x \in B\}.$
 
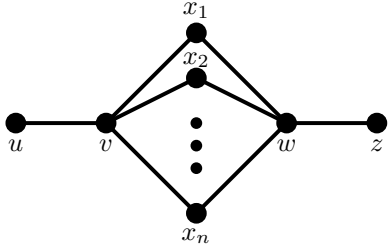
\begin{figure}
    \centering
\begin{tikzpicture}[scale=.6]
\draw[fill=black!100,line width=1] (0,-2) circle (.2);
\draw[fill=black!100,line width=1] (0,2) circle (.2);
\draw[fill=black!100,line width=1] (2,0) circle (.2);
\draw[fill=black!100,line width=1] (-2,0) circle (.2);
\draw[fill=black!100,line width=1] (-4,0) circle (.2);
\draw[fill=black!100,line width=1] (4,0) circle (.2);
\draw[fill=black!100,line width=1] (0,1) circle (.2);
\draw[fill=black!100,line width=1] (0,0) circle (.1);
\draw[fill=black!100,line width=1] (0,-.5) circle (.1);
\draw[fill=black!100,line width=1] (0,-1) circle (.1);
\draw[black!100,line width=1.5pt] (0,-2) -- (2,0) -- (0,2) -- (-2,0) -- (0,-2);
\draw[black!100,line width=1.5pt] (-2,0) -- (0,1) -- (2,0);
\draw[black!100,line width=1.5pt] (-4,0) -- (-2,0);
\draw[black!100,line width=1.5pt] (4,0) -- (2,0);
\draw (-4,-.5) node{$u$};
\draw (-2,-.5) node{$v$};
\draw (2,-.5) node{$w$};
\draw (4,-.5) node{$z$};
\draw (0,1.4) node{$x_2$};
\draw (0,2.5) node{$x_1$};
\draw (0,-2.5) node{$x_n$};
\end{tikzpicture}
    \caption{A graph $G$ for which no minimum twin cover is determining~\cite{BCKLPR2020b}.}
    \label{fig:C4Pend}
\end{figure}

We call a minimum size subset of $V(G)$ containing at least one vertex from every pair of twin vertices a {\em minimum twin cover}. If $G$ has no twins, then a minimum twin cover is the empty set. Every determining set of a graph $G$ must contain all but one vertex in every collection of mutual twins, i.e.\ a minimum twin cover. Thus, if a minimum twin cover is a determining set, it must be a minimum size determining set. In this case, every minimum twin cover is a minimum size determining set, since all minimum twin covers have the same cardinality. To see that this is true, note that a minimum twin cover contains all but one vertex from every equivalence class $[x]$, $x \in V(G)$. In other words, every minimum twin cover has cardinality $|V(G)|-|V(\tG)|$. Figure~\ref{fig:C4Pend} shows an example of a graph in which minimum twin covers are not determining sets.

The following three results, from~\cite{BCKLPR2020b}, establish relationships between minimum twin covers of $G$, minimum determining sets for $G$, and minimum determining sets for $\tG$.   

\begin{cor}\label{cor:TtotT} Let $G$ be a graph.  If $S$ is a determining set for $G$, then $\tS$ is a determining set for $\tG$.\end{cor}

Observe that a superset of a determining set is still a determining set. Therefore, if  $\widetilde D$ is a determining set for $\tG$, then  $\tT \cup \widetilde D$ is a determining set for $\tG$ containing $\tT$. Similarly, if $S$ is a determining set for $G$ containing $T$, then $\tS$ is a determining set for $\tG$. 

\begin{thm} \label{thm:DettG}~{\rm \cite{BCKLPR2020b}}
Let $T$ be a minimum twin cover of $G$ and let $\tS$ be a determining set for $\tG$ containing $\tT$. Let $R=\{x\in V(G) \ | \ [x]\in \tS\setminus \tT\}$. Then $S=T \cup R$ is a determining set for $G$. Furthermore, if $\tS$ is of minimum size among determining sets for $\tG$ that contain $\tT$, then $S$ is a minimum size determining set for $G$.

\end{thm}

In particular, if $\tS = \tT$, then $T$ is a minimum size determining set for $G$. This observation allows us to see that Theorem~\ref{thm:DettG} yields natural bounds on $\det(G)$ in terms of $|T|$ and $\det(\tG)$. 

\begin{cor}\label{cor:twinBounds}
Let $T$ be a minimum twin cover of $G$. Then 
\[
|T| \leq \det(G) \leq |T| {+} \det(\tG),
\] with both bounds sharp.
\end{cor}

To use these results in the context of generalized Mycielskian constructions, we now investigate how applying the generalized Mycielskian construction affects the size of a minimum twin cover as well as the relationship between $\det(\tG)$ and $\det(\widetilde \gm)$. This is done in the next two lemmas.

In~\cite{BCKLPR2020b}, it is shown that if $T$ is a minimum twin cover in a graph $G$ with no isolated vertices, the set consisting of vertices in $T$ and all of their shadows, 
\[
\Tt =  \{u_i^s \mid v_i \in T, \, 0 \le s \le t \},
\]
is a minimum twin cover of $\gm$ of size $(t{+}1)|T|$. 

If $G$ has isolated vertices, then an extra $t-1$ vertices need to be added to $\Tt$ to obtain a minimum twin cover of $\gm$. This is because if $G$ has $\ell \geq 1$ isolated vertices, only $\ell - 1$ of them are in any minimum twin cover of $G$. Thus, there is always one isolated vertex in $G$ not in $t$. If $u$ is the isolated vertex in $G$ not in the minimum twin cover, $\gm$ has $t-1$ shadows of $u$ that are mutually twin with the other isolated vertices. This behavior is different from non-isolated vertices, which cannot have twins that are on different levels.

\begin{lemma}\label{lem:twinrepsmuIso}
Let $T$ be a minimum twin cover of $G$.
If $G$ has at least one isolated vertex and $u$ is the isolated vertex that is not in $T$, then for $t \geq 1$,
\[ \Tt\cup \{u^s \mid  1\le s \le t{-}1\} = \{ u_i^s \mid v_i \in T, 0 \le s \le t\} \cup \{u^s \mid 1\le s \le t{-}1\} \]
is a minimum twin cover of $\gm$ of size $(t{+}1)|T| {+} t{-}1$. 
\end{lemma}

\begin{proof} First suppose  $G = \ell K_1$ for $\ell \ge 1$. Then $\gm =K_{1,\ell} + t\ell K_1$. Since $T$ must contain $\ell{-}1$ of the isolated vertices, the set $\Tt = \{ u_i^s \mid v_i \in T, 0 \le s \le t\}$ contains $t(\ell{-}1)$ of the isolated vertices and $\ell{-}1$ of the leaves of the $K_{1,\ell}$. The set $\{u^s \mid 1\le s \le t{-}1\}$ contains an additional $t{-}1$ of the isolated vertices in $\gm$. Thus, their union gives a minimum twin cover of $\gm$. 

Next, suppose $G = H + \ell K_1$, where $H$ has no isolated vertices and $\ell \ge 1$. Then $\gm = C {+} t\ell K_1$ where $C$ is $\mu_t(H)$ with an additional $\ell$ pendant vertices adjacent to $w$. By Observation~\ref{obs:twins} and  the fact that the isolated vertices are  mutually twin in $G$, every pair of twins in $C$ has at least one member in $\Tt$. Furthermore, since $T$ excludes exactly one vertex from each set of mutual twins in $G$, $\Tt\cap V(C)$ is a minimum twin cover of $C$. 

Finally, a minimum twin cover of $\gm$ must contain $t\ell-1$ of the isolated vertices. The sets $\Tt\setminus V(C)$ and $\{u^s | 1\le s \le t{-}1\}$ give these.\end{proof}

If a graph $G$ has twins but no isolated vertices, then it is shown in~\cite{BCKLPR2020b} that the process of applying the canonical quotient map commutes with the process of applying the generalized Mycielskian construction; that is,  \[ \mu_t(\tG) = \widetilde{\gm}. \]
However, this is not the case for graphs with isolated vertices. As before, an extra $t-1$ vertices are needed.
 
\begin{lemma}\label{lem:commutesIso}
 Let $G$ be a graph with isolated vertices such that $G \ne \tG$. Then for $t \geq 1$, \[\mu_t(\tG) = \widetilde{\gm}{+}(t{-}1)K_1.\]
\end{lemma}

\begin{proof}

Let $G = H + \ell K_1$, where $H$ has no isolated vertices and $\ell \ge 1$. Then $\gm = C + t\ell K_1$, where $C$ is a connected graph consisting of a copy of $\gh$ with $\ell$ pendant vertices adjacent to $w$. Thus $\widetilde{\gm}= \widetilde C + K_1$. 
It is straightforward to see that $\widetilde C$ consists of a copy of $\widetilde{\gh} = \mu_t(\tH)$ with a single extra pendant vertex adjacent to $[w]$. 

On the other hand, $\tG = \tH + K_1$, where $\tH$ has no isolated vertices. So $\mu_t(\tG)$ is a connected graph consisting of a copy of $\mu_t(\tH)$ with a single pendant vertex adjacent to $[w]$, which is $\widetilde C$,  plus $t$ isolated vertices. 
Thus, $\mu_t(\tG) = \widetilde{\gm}{+}(t{-}1)K_1$.
\end{proof}

 In~\cite{BCKLPR2020b}, it is shown that if a graph $G$ has twins, then we can find the determining number of $\gm$ by using a determining set in $\tG$ that contains $\tT$. The result is the following.

\begin{thm}~{\rm \cite{BCKLPR2020b}}\label{thm:twinDet}
 Let $G$ be a graph with twins and no isolated vertices. Let $T$ be a minimum twin cover of $G$. Then for $t \geq 1$,

 \[\det(\gm)=t|T| + \det(G).\]

\end{thm}

If $G$ has isolated vertices, then by Lemma~\ref{lem:twinrepsmuIso}, a minimum twin cover of $\gm$ requires an extra $t-1$ vertices beyond the vertices in a minimum twin cover of $G$ and their shadows, and so it is not surprising that its determining number must also increase by $t-1$.
 
 \begin{thm}\label{thm:twinswithiso}
Let $G$ be a graph with twins and isolated vertices and let $T$ be a minimum twin cover of $G$. Then for $t\geq 1$,

\[ \det(\gm) = t|T| + \det(G) + t - 1.\]
\end{thm}

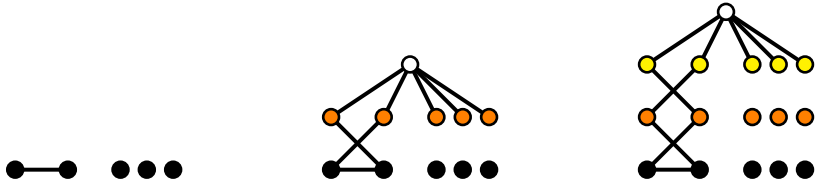
\begin{figure}[h]
    \centering
\begin{tikzpicture}[scale=.7]

\draw[black!100,line width=1.5pt] (0,0) -- (1,0); 
\draw[fill=black!100,line width=1] (0,0) circle (.15);
\draw[fill=black!100,line width=1] (1,0) circle (.15);
\draw[fill=black!100,line width=1] (2,0) circle (.15);
\draw[fill=black!100,line width=1] (2.5,0) circle (.15);
\draw[fill=black!100,line width=1] (3,0) circle (.15);
\begin{scope}[shift={(6,0)}]
\draw[black!100,line width=1.5pt] (0,0) -- (1,0); 
\draw[fill=black!100,line width=1] (0,0) circle (.15);
\draw[fill=black!100,line width=1] (1,0) circle (.15);
\draw[fill=black!100,line width=1] (2,0) circle (.15);
\draw[fill=black!100,line width=1] (2.5,0) circle (.15);
\draw[fill=black!100,line width=1] (3,0) circle (.15);
\draw[black!100,line width=1.5pt] (0,0) -- (1,1); 
\draw[black!100,line width=1.5pt] (0,1) -- (1,0); 
\draw[black!100,line width=1.5pt] (1.5,2) -- (0,1); 
\draw[black!100,line width=1.5pt] (1.5,2) -- (1,1); 
\draw[black!100,line width=1.5pt] (1.5,2) -- (2,1); 
\draw[black!100,line width=1.5pt] (1.5,2) -- (2.5,1); 
\draw[black!100,line width=1.5pt] (1.5,2) -- (3,1); 
\draw[fill=orange!100,line width=1] (0,1) circle (.15);
\draw[fill=orange!100,line width=1] (1,1) circle (.15);
\draw[fill=orange!100,line width=1] (2,1) circle (.15);
\draw[fill=orange!100,line width=1] (2.5,1) circle (.15);
\draw[fill=orange!100,line width=1] (3,1) circle (.15);
\draw[fill=white!100,line width=1] (1.5,2) circle (.15);
\end{scope}
\begin{scope}[shift={(12,0)}]
\draw[black!100,line width=1.5pt] (0,0) -- (1,0); 
\draw[fill=black!100,line width=1] (0,0) circle (.15);
\draw[fill=black!100,line width=1] (1,0) circle (.15);
\draw[fill=black!100,line width=1] (2,0) circle (.15);
\draw[fill=black!100,line width=1] (2.5,0) circle (.15);
\draw[fill=black!100,line width=1] (3,0) circle (.15);
\draw[black!100,line width=1.5pt] (0,0) -- (1,1); 
\draw[black!100,line width=1.5pt] (0,1) -- (1,0);
\draw[black!100,line width=1.5pt] (0,1) -- (1,2); 
\draw[black!100,line width=1.5pt] (0,2) -- (1,1);
\draw[black!100,line width=1.5pt] (1.5,3) -- (0,2); 
\draw[black!100,line width=1.5pt] (1.5,3) -- (1,2); 
\draw[black!100,line width=1.5pt] (1.5,3) -- (2,2); 
\draw[black!100,line width=1.5pt] (1.5,3) -- (2.5,2); 
\draw[black!100,line width=1.5pt] (1.5,3) -- (3,2); 
\draw[fill=orange!100,line width=1] (0,1) circle (.15);
\draw[fill=orange!100,line width=1] (1,1) circle (.15);
\draw[fill=orange!100,line width=1] (2,1) circle (.15);
\draw[fill=orange!100,line width=1] (2.5,1) circle (.15);
\draw[fill=orange!100,line width=1] (3,1) circle (.15);
\draw[fill=yellow!100,line width=1] (0,2) circle (.15);
\draw[fill=yellow!100,line width=1] (1,2) circle (.15);
\draw[fill=yellow!100,line width=1] (2,2) circle (.15);
\draw[fill=yellow!100,line width=1] (2.5,2) circle (.15);
\draw[fill=yellow!100,line width=1] (3,2) circle (.15);
\draw[fill=white!100,line width=1] (1.5,3) circle (.15);
\end{scope}
\end{tikzpicture}
    \caption{The graphs $G=K_2+3K_1$, $\mu(G)$ and $\mu_2(G)$. }
    \label{fig:K2tK1ex}
\end{figure}

\begin{proof}
First suppose that $G = \ell K_1$ with $\ell \geq 2$ since $G \neq \tG$. Then $T$ contains $\ell-1$ of the isolated vertices of $G$; such a set is also a minimum determining set for $G$ and so $\det(G) = |T|$.

Observe that $\gm = K_{1,\ell} + t\ell K_1$. It is clear that a minimum determining set for $\gm$ consists of $\ell-1$ of the pendant vertices of  $K_{1,\ell}$, together with $t \ell - 1$ of the isolated vertices. It follows that \[\det(\gm) = (t+1)\ell - 2  = (t+1)|T| +t-1 = t|T| + \det(G) + t - 1.\]

Now suppose $G  = H+\ell K_1$ with $\ell \ge 1$ and where $H$ has at least one edge and no isolated vertices. By construction, $\gm = C+ t\ell K_1$, where $C$ consists of $\mu_t(H)$ with $\ell$ pendant vertices adjacent to $w$.  Then, \[\det(\gm) = \det(C) + \det(t\ell K_1) = \det(C) + t\ell - 1.\] 

Since $H$ has no isolated vertices, a minimum twin cover $T$ of $G$ must consist of a minimum twin cover $T_H$ of $H$ plus $\ell - 1$ of the isolated  vertices. Thus, $|T| = |T_H| + \ell-1$.

In the following we show that $\det(C) = t|T_H| + \det(H) + \ell - 1$. We divide into two cases, depending on whether $H$ is a star graph of the form $K_{1,m}$ for $m \ge 1.$ 

First suppose $H \neq K_{1,m}$. Since $G = H + \ell K_1$, $\tG = \tH + K_1$. Let $\tA$ be a minimum determining set for $\tG$ containing $\tT$ and let $R =  \{x \in V(G) \mid [x]\in \tA \setminus \tT\}$. Notice that, if $\ell = 1$, one isolated vertex is not in $T$ and its equivalence class is not in $\tA$ and if $\ell \geq 2$, the isolated vertex in $\tG$ is in $\tT$. Thus, $R \subseteq V(H)$.

Since $\tA$ is a minimum size determining set for $\tG$ containing $\tT$, $\tA \cap V(\tH) = R \cup \widetilde{T_H}$
 is a minimum size determining set for $\tH$ containing $\widetilde{T_H}$. By the proof of Theorem~\ref{thm:twinDet}, we then get that $R \cup T_H$ is a determining set for $\gh$ and thus, $\det(\gh) = (t+1)|T_H| + |R| = t|T_H| + \det(H)$.
 
Now, since $C$ is a copy of $\gh$ with $\ell$ pendant vertices adjacent to $w$, any minimum determining set for $C$ must include $\ell - 1$ of the pendant vertices. It follows that $\det(C) \leq \det(\gh) + \ell - 1$. Since $H \neq K_{1,m}$, any automorphism of $\gh$ fixes $w$ by Lemma~\ref{lem:phiwisv}. Thus, any minimum determining set for $\gh$ does not contain $w$. Thus, by Lemma~\ref{lem:detghvsC}, \[\det(C) = \det(\gh) + \ell - 1 = t|T_H| + \det(H) + \ell - 1.\]

Now suppose that $H$ is a star graph.  We begin with the case that $H = K_{1,m}$ for $m \geq 2$. Then $\tH = K_2$. A minimum twin cover $T$ of $G$ consists of $m-1$ of the pendant vertices in $H$ and $\ell-1$ of the isolated vertices. In particular, $T_H = T \cap V(H)$ is a minimum twin cover of $H$ that is also a determining set for $H$. Thus, $\det(H) = |T_H|$.  

 Since $H$ has no isolated vertices, $(T_H)_t$ is a minimum twin cover of $\mu_t(H)$.
 We next show that $(T_H)_t$ is also a determining set for  $\mu_t(H)$.
 
 In $\mu_t(H)$ the minimum twin cover contains $m-1$ vertices on each level. Let $v$ be the vertex in $H$ of degree $m$, and $v_1,\dots,v_{m-1}$ be the vertices in $T_H$. In the following, we show that any automorphism of $\mu_t(H)$ that fixes $(T_H)_t$ fixes $w$. Observe that the neighbors of $w$ in $\mu_t(H)$ that are in $(T_H)_t$ are $u_1^t,\dots,u_{m-1}^t$. The only other vertex with this property is $u^{t-1}$. However, the vertex $u^{t-1}$ has degree $2m$ and $w$ has degree $m+1$. Therefore, any automorphism of $\mu_t(H)$ that fixes $(T_H)_t$ must also fix $w$. Since automorphisms preserve distance, this implies that such an automorphism preserves levels. Within a level, fixing the vertices in $(T_H)_t$ fixes everything. Therefore, $(T_H)_t$ is a determining set for $\mu_t(H)$.
 
Since a minimum twin cover is a determining set,  all minimum size determining sets are minimum twin covers, each of which excludes $w$. 
 In particular, $(T_H)_t$ is a minimum size subset of $V(\gh) \setminus \{w\}$ such that $(T_H)^{(t)} \cup \{w \}$ is a determining set.
 By Lemma~\ref{lem:detghvsC},  $(T_H)_t$ plus  $\ell - 1$ of the pendant vertices adjacent to $w$ constitute a minimum size determining set for $C$. 
  Thus, $\det(C) = (T_H)_t + \ell - 1 = (t+1)(T_H) + \ell - 1$, as desired.
 
  Finally, we consider the case where $H$ is the star graph $K_{1,1} = K_2.$ Since we are assuming that $G$ has twins,  $G = K_2 + \ell K_1$ for some $\ell \ge 2$. 
A minimum twin cover $T$ of $G$ consists of $\ell - 1$ of the isolated vertices, but no vertices of $H$. Thus, $T_H$ is empty and so $|T_H| = 0$. Also, $T$ is not a determining set for $G$; instead the determining set for $G$ consists of one of the two vertices in $K_2$ together with $T$. It follows then that $\det(H) = 1$.

In this case, $\gm = C+  t\ell K_1$, 
where $C$ is a copy of $\mu_t(K_2) = C_{2t+3}$ with an additional $\ell$ pendant vertices $x_1, \dots, x_\ell$ adjacent to $w$. This is illustrated in Figure~\ref{fig:K2tK1ex}. If $y \ne w$ is any vertex in the cycle, then $\{y\} \cup \{w\}$ is a minimum size determining set for $\mu_t(K_2)$ and so by Lemma~\ref{lem:detghvsC}, $\{y, x_2, \dots , x_\ell\}$ is a minimum size determining set for $C$.  Thus, 
\[
\det(C) = \ell = t |T_H| + \det(H) + (\ell -1).
\]

Finally,  
\begin{align*}
\det(\gm) &= \det(C) + t\ell-1 \\
&= \big [t |T_H| + \det(H) + (\ell -1) \big ]+ t\ell - 1\\
&= \big [t |T_H| + t\ell  - t \big ]  +\big [ \det(H) + (\ell-1) \big ] + t - 1\\
&= \big [t |T_H| + t(\ell  - 1) \big ] + \big [ \det(H) + (\ell-1) \big ] + t - 1\\
&= t|T| + \det(G) + t-1,
\end{align*}
as desired.
 \end{proof}

 If a minimum twin cover $T$ is a determining set for $G$, then by Corollary~\ref{cor:TtotT}, $\tT$ is a determining set for $\tG$. Thus $|T| =\det(G)$ and $|R|= 0$, and so Theorem~\ref{thm:twinswithiso}  implies the following corollary. 

 \begin{cor}\label{cor:twincovdet}
Let $G$ be a graph with isolated vertices such that $G \neq \tG$. If $G$ has a minimum twin cover that is also a determining set, then for $t \ge 1$,
    \[\det(\gm) = (t{+}1)\det(G)+t-1.\]

\end{cor}

Note that for $R$ in the hypothesis of Theorem~\ref{thm:twinswithiso}, $|R| \le \det(\tG)$. If no minimum determining set for $\tG$ contains any vertex of $\tT$, then $|R| = \det(\tG)$. This is the case, for example, for the graph $G$ in Figure~\ref{fig:C4Pend} or for $G = K_2 + \ell K_1$ with $\ell \ge 2$. Combining this observation with Corollary~\ref{cor:twincovdet} gives the following corollary, bounding $\det(\gm)$ for all graphs $G \neq \tG$ with isolated vertices. 

\begin{cor}\label{cor:twinDetboundIso}
Let $G$ be a graph with isolated vertices such that $G \neq \tG$.
Let $T$ be a minimum twin cover of $G$. 
Then for $t \geq 1$, \[
(t{+}1)|T|{+}t{-}1 \leq \det(\gm) \leq 
\det(\tG) {+}(t{+}1)|T| {+}t{-}1.
\]
Both bounds are sharp.
\end{cor}

 We end our discussion on the determining number of the generalized Mycielskians of graphs with isolated vertices by combining 
Theorem~\ref{thm:twinswithiso} and Theorem~\ref{twinfreeDetIso}. Observe that if $G$ is twin-free, then the minimum twin cover is empty and $|R| = \det(G)$.

\begin{thm}\label{thm:isolatescombo}
Let $G$ be a graph with isolated vertices and let $T$ be a (possibly empty) minimum twin cover of $G$. Among all determining sets for $\tG$ containing $\tT$, let $\tA$ be one of minimum size. Let $R =  \{x \in V(G) \mid [x]\in \tA \setminus \tT\}$. Let $t\geq 1$.

\begin{enumerate}[(i)]
   \item If $G = K_1$ then, $\det(G) = 0$ and $\det(\gm) = t.$
   \item If $G \neq K_1$, then $\det(\gm) =  t|T| + \det(G) + t-1$.
   
   \end{enumerate}

\end{thm}

\section{Open Problems}\label{sec:Open}

There are a significant number of results on distinguishing graph families in the literature. However, relatively little work has been done on the cost of 2-distinguishing or on finding determining numbers. 

Some of the existing work on determining numbers is related to Cartesian products. The \emph{Cartesian product} of graphs $G$ and $H$ has vertex set $V(G)\times V(H)$ with an edge between vertices $(x,u)$ and $(y,u)$ if $x$ is adjacent to $y$ in $G$ and between vertices $(x,u)$ and $(x,v)$ if $u$ is adjacent to $v$ in $H$. The \emph{Cartesian power} $G^k$ of $G$ is the Cartesian product of $G$ with itself $k$ times. A graph $G$ is \emph{prime} with respect to the Cartesian product if it cannot be written as the Cartesian product of two smaller graphs. 

\begin{thm}
\cite{Bo2013a} If $G^k$ is a 2-distinguishable Cartesian power of a prime connected graph $G$ on at least three vertices with $\det(G) \le k$ and $\max\{2,\det(G)\}<\det(G^k)$, then $\rho(G^k)\in\{\det(G^k),\det(G^k) + 1\}$.
\end{thm}

\begin{problem}
\cite{Bo2013a} Classify the 2-distinguishable Cartesian powers $G^k$ so that $\rho(G^k) = \det(G^k)$, or so that $\rho(G^k)=\det(G^k)+1$. 
\end{problem}

\begin{problem}
\cite{Bo2013a} 
More generally, classify the 2-distinguishable graphs with \[\rho(G) \in \{\det(G), \det(G)+1\}.\] 
\end{problem}

One could choose to focus on either $\rho(G) = \det(G)$ or $\rho(G) = \det(G) +1$ for the classification instead. 

\smallskip 

 For Mycielskian graphs, some work has been done on the edge version of distinguishing. A coloring of the edges of a graph $G$ with colors from $\{1,\ldots, d\}$ is called a \emph{$d$-distinguishing edge coloring} if no nontrivial automorphism of $G$ preserves the color classes. No such coloring exists when $G$ has $K_2$ as a component nor when $G$ has two or more $K_1$ components. For all other graphs $G$, the \emph{distinguishing index} of $G$, denoted $\dist'(G)$, is the smallest number of colors in a $d$-distinguishing edge coloring of $G$.
 
 \begin{thm}
 \cite{AlSo2019b} Let $G$ be a twin free graph with at least three vertices, no $K_2$ component and at most one $K_1$ component. Then $\dist'(\mug) \le \dist'(G)+1$. 
 \end{thm}
 
 \begin{conj}
  \cite{AlSo2019b} Let $G$ be a connected graph with at least three vertices. Then $\dist'(\mug) \le \dist'(G)$ except for a finite number of graphs. 
 \end{conj}
 
 \begin{problem}
 Express $\dist'(\mu_t(G))$ in terms of $\dist'(G)$ and possibly other parameters of $G$. 
 \end{problem}
 
 \begin{problem}
For each graph $G$, determine $t$ for which $\dist'(\mu_t(G))=2$.  
 \end{problem}

\section{Acknowledgments}

This work is a result of a collaboration made possible by the Institute for Mathematics and its Applications' Workshop for Women in Graph Theory and Applications, August  2019.

\FloatBarrier
\bibliographystyle{plain}
\bibliography{Chapter}

\end{document}